\documentclass[11pt]{article}

\usepackage[a4paper,margin=1in]{geometry}
\usepackage{amsmath,amssymb,amsthm,mathtools}
\usepackage{mathrsfs}
\usepackage{enumitem}
\usepackage{hyperref}
\usepackage{cleveref}
\usepackage{microtype}

\newcommand{\F}{\mathcal{F}}

\newcommand{\floor}[1]{\left\lfloor #1 \right\rfloor}

\title{New approaches to the
Erd\H{o}s Distinct Subset Sums Problem}

\date{}

\author{
Aren Hasratyan\\
\small Sorbonne University Abu Dhabi, Abu Dhabi, UAE\\
\small \texttt{arenhasratyan2006@gmail.com}
\and
Mikl\'os Ruszink\'o\\
\small Sorbonne University Abu Dhabi, SAFIR, CEFAMSDS, Abu Dhabi, UAE\\
\small Alfr\'ed R\'enyi Institute of Mathematics, Budapest, Hungary\\
\small \texttt{ruszinko@renyi.hu}
}

\begin{document}

\maketitle


\begin{abstract}
    We present a new method of reproving the best known lower bound of the Erd\H{o}s distinct subset sums problem. Although our lower bound is asymptotically the same as in \cite{DubroffFoxXu2021}, it is slightly better for the case when \(n\) is odd. Besides, we also provide a new lower bound on the second largest element \(a_{n-1}\).
\end{abstract}

\section{Introduction}
Let
\[
A=\{a_1<a_2<\cdots<a_n\}
\]
be a set of positive integers whose subset sums are all distinct.
The Erd\H{o}s distinct subset sums problem asks how small the
largest element \(a_n\) can be. In particular, Erd\H{o}s asked
whether there exists an absolute constant \(c>0\) such that
\(a_n\geq c2^n\) and offered \(500\$\) for proving or disproving \cite{Erdos1956}, which was claimed to be disproved by GPT-6 Astra in \(2026\). The classical variance method
gives a lower bound of order \(2^n/\sqrt n\), with subsequent
improvements by Elkies \cite{Elkies1986}, Bae \cite{Bae1996},
and Aliev \cite{Aliev2008}. Dubroff, Fox, and Xu obtained the
asymptotic lower bound
\[
a_n\geq
\left(\sqrt{\frac{2}{\pi}}-o(1)\right)
\frac{2^n}{\sqrt n}
\]
and provided two short proofs of this estimate
\cite{DubroffFoxXu2021}. In the opposite direction, before GPT-6 Astra, the best known upper bound was by Bohman, who
constructed sets with
\[
a_n<0.22002\,2^n
\]
for all sufficiently large \(n\) \cite{Bohman1998}.

\section{Main results}
Let
\[
A=\{a_1<a_2<\cdots<a_n\}
\]
be a set of positive integers with distinct subset sums, and put
\[
T=\sum_{j=1}^{n-1}a_j.
\]
Write the subset sums of \(a_1,\ldots,a_{n-1}\) in increasing
order:
\[
0=t_1<t_2<\cdots<t_{2^{n-1}}=T.
\]

\noindent\textbf{Theorem 1.}
If $a_1<a_2<\cdots<a_n$ are positive integers with distinct subset sums, then 

\[
a_n \geq 2\binom{n-1}{\left\lfloor \frac{n-1}{2}\right\rfloor}.
\]

In our proof, we combine ideas from \cite{DubroffFoxXu2021} with the following observations. If $a_n$ is $\Omega (2^n)$, then the theorem is proven, i.e., we may assume that $a_n = o(2^n)$.

\medskip
\noindent\textbf{Lemma 1. }
If \(i+a_n\leq 2^{n-1}\), then
\[
t_{i+a_n}-t_i\geq 2a_n.
\]

\begin{proof}
Consider the \(a_n+1\) integers
\[
t_i,t_{i+1},\ldots,t_{i+a_n}.
\]
By the pigeonhole principle, two of them, say \(t_k\) and
\(t_\ell\), with \(k<\ell\), are congruent modulo \(a_n\).
Hence
\[
t_\ell-t_k=ma_n
\]
for some positive integer \(m\).

If \(m=1\), then
\[
t_\ell=t_k+a_n.
\]
The left-hand side is a subset sum of
\(a_1,\ldots,a_{n-1}\), while the right-hand side is a subset
sum involving \(a_n\). This contradicts the distinctness of all
subset sums. Therefore \(m\geq2\), and
\[
t_{i+a_n}-t_i
\geq t_\ell-t_k
\geq2a_n.
\]
\end{proof}

\medskip
\noindent\textbf{Lemma 2.}
Define
\[
\F=
\left\{
I\subseteq[n-1]:
\sum_{j\in I}a_j<\frac{T}{2}
\right\}=
\left\{
t_i:
t_i<\frac{T}{2}
\right\}.
\]
Then
\[
|\F|=2^{n-2}.
\]

\begin{proof}
For any \(1\leq i \leq 2^{n-1}\), \(t_i + t_{2^{n-1}+1-i}=T\). Since all subset sums are distinct, no subset sum equals \(T/2\). So, exactly one member of each
pair belongs to \(\F\). Therefore,
\[
|\F|=2^{n-2}.
\]
\end{proof}

Now, like in \cite{DubroffFoxXu2021}, identify subsets of \([n-1]\) with vertices of the
\((n-1)\)-dimensional Boolean cube. For a family
\(\F\subseteq2^{[n-1]}\), let \(\partial\F\) denote its
external vertex boundary. Harper's vertex-isoperimetric
inequality implies that if
\[
|\F|=2^{n-2},
\]
then
\[
|\partial\F|
\geq
\binom{n-1}{\floor{(n-1)/2}}.
\]

From this, we can conclude that 
\[
\#\left\{
t_i:
\frac{T}{2}-t_i<a_{n-1}
\right\}
\geq
\binom{n-1}{\floor{(n-1)/2}}.
\] 

As \(t_i\)-s are symmetric with respect to \(\frac{T}{2}\) and \(a_n>a_{n-1}\), we can conclude 

\[
\#\left\{
t_i:
\left|t_i-\frac{T}{2}\right|<a_n
\right\}
\geq
\#\left\{
t_i:
\left|t_i-\frac{T}{2}\right|<a_{n-1}
\right\}
\geq
2\binom{n-1}{\floor{(n-1)/2}}.
\]

\medskip
\noindent\textbf{Lemma 3.}

\[
a_n
\geq
\#\left\{
t_i:
\left|t_i-\frac{T}{2}\right|<a_n
\right\}.
\]

\begin{proof}
Consider the interval
\[
J=\left(\frac{T}{2}-a_n,\,
         \frac{T}{2}+a_n\right),
\]
whose length is \(2a_n\).

If \(J\) contained \(a_n+1\) consecutive subset sums
\[
t_i,t_{i+1},\ldots,t_{i+a_n},
\]
then
\[
t_{i+a_n}-t_i<2a_n,
\]
contradicting the first lemma. Therefore \(J\) contains at most
\(a_n\) subset sums.

\end{proof}

Hence, we can conclude from those lemmas that 
\[
a_n
\geq
2\binom{n-1}{\floor{(n-1)/2}},
\]
which is exactly the statement of the theorem.

\medskip
\noindent\textbf{Asymptotic consequence.}

Using the central-binomial-coefficient estimate
\[
\binom{m}{\floor{m/2}}
=
\left(\sqrt{\frac{2}{\pi}}-o(1)\right)
\frac{2^m}{\sqrt m},
\]
the conditional bound in the theorem yields
\[
a_n\geq
\left(\sqrt{\frac{2}{\pi}}-o(1)\right)
\frac{2^n}{\sqrt {n-1}}\geq
\left(\sqrt{\frac{2}{\pi}}-o(1)\right)
\frac{2^n}{\sqrt n}.
\]
This is of the same asymptotic order as the established lower
bound discussed by Dubroff, Fox, and Xu
\cite{DubroffFoxXu2021}. It is also worth mentioning that this lower bound is better for odd \(n\), since if \(n\) is odd, then \( 2\binom{n-1}{\floor{(n-1)/2}} >
\binom{n}{\floor{n/2}}
\), as \( 2\binom{n-1}{\floor{(n-1)/2}} -
\binom{n}{\floor{n/2}}=C_{\frac{n-1}{2}}
\), where \(C_{\frac{n-1}{2}}\) is the \(\frac{n-1}{2}\)-th Catalan number.


\section{Some remarks on the second largest element}

For this subsection, again distinguish the subset sums of all \(n\)
elements from those of the first \(n-1\) elements. Write
\[
0=t_1<t_2<\cdots<t_{2^{n-1}}=T
=\sum_{i=1}^{n-1} a_i.
\]

\noindent\textbf{Theorem 2.}

If $a_1<a_2<\cdots<a_n$ are positive integers with distinct subset sums, then 

\[
{
a_{n-1}\geq
\left(\frac{2}{\sqrt{3}}-o(1)\right)\frac{2^{n-1}}{\sqrt{n-1}}
}.
\]

\medskip
\noindent\textbf{Lemma 4.}
\[
\frac{2^{2(n-1)}}{4}\sum_{i=1}^{n-1} a_i^2
= \sum_{i>j}(t_i-t_j)^2.
\]

\begin{proof}
Since \(t_1,\ldots,t_{2^{n-1}}\) are all the subset sums of
\(a_1,\ldots,a_{n-1}\), we have
\[
\begin{aligned}
\sum_{i>j}(t_i-t_j)^2
&=\frac12\sum_{i=1}^{2^{n-1}}
\sum_{j=1}^{2^{n-1}}(t_i-t_j)^2\\
&=\frac12\sum_{i=1}^{2^{n-1}}
\sum_{j=1}^{2^{n-1}}(t_i^2+t_j^2-2t_it_j)\\
&=2^{n-1}\sum_{i=1}^{2^{n-1}}t_i^2
-\left(\sum_{i=1}^{2^{n-1}}t_i\right)^2.
\end{aligned}
\]

Each \(a_i\) appears in exactly \(2^{n-2}\) subset sums. Hence
\[
\sum_{i=1}^{2^{n-1}}t_i=2^{n-2}\sum_{i=1}^{n-1}a_i
=2^{n-2}T.
\]

Moreover, expanding the squares of all subset sums, each
\(a_i^2\) appears \(2^{n-2}\) times, while each product
\(a_ia_j\), with \(i<j\), appears \(2^{n-3}\) times.
Therefore,
\[
\begin{aligned}
\sum_{i=1}^{2^{n-1}}t_i^2
&=2^{n-2}\sum_{i=1}^{n-1}a_i^2
+2^{n-2}\sum_{1\le i<j\le n-1}a_ia_j\\
&=2^{n-3}\left(
\sum_{i=1}^{n-1}a_i^2+T^2
\right),
\end{aligned}
\]
where we used
\[
T^2=\sum_{i=1}^{n-1}a_i^2
+2\sum_{1\le i<j\le n-1}a_ia_j.
\]

Substituting these expressions gives
\[
\begin{aligned}
\sum_{i>j}(t_i-t_j)^2
&=2^{n-1}\cdot 2^{n-3}
\left(\sum_{i=1}^{n-1}a_i^2+T^2\right)
-2^{2n-4}T^2\\
&=2^{2n-4}\sum_{i=1}^{n-1}a_i^2\\
&=\frac{2^{2(n-1)}}{4}
\sum_{i=1}^{n-1}a_i^2.
\end{aligned}
\]
\end{proof}

From Lemma \(1\), we have
\[
t_i-t_j \geq 2a_n
\left\lfloor \frac{i-j}{a_n}\right\rfloor.
\]

Therefore,
\[
\begin{aligned}
\frac{2^{2(n-1)}}{4}\sum_{i=1}^{n-1}a_i^2
&\geq 4a_n^2
\sum_{i>j}
\left\lfloor\frac{i-j}{a_n}\right\rfloor^2\\
&\geq \left(1-o(1)\right) 4a_n^2
\sum_{k=1}^{\lfloor2^{n-1}/a_n\rfloor}
k^2a_n(2^{n-1}-ka_n)\\
&=\left(1-o(1)\right)4a_n^3
\sum_{k=1}^{\lfloor 2^{n-1}/a_n \rfloor}
\left(2^{n-1}k^2-a_nk^3\right)\\
&\geq \left(1-o(1)\right)4a_n^3\left(
\frac{2^{n-1}}{3}
\frac{2^{3(n-1)}}{a_n^3}
-\frac{a_n}{4}
\frac{2^{4(n-1)}}{a_n^4}
\right)\\
&=\left(1-o(1)\right)\frac{2^{4(n-1)}}{3}.
\end{aligned}
\]

Hence,
\[
\sum_{i=1}^{n-1}a_i^2
\geq (1-o(1)) \frac{2^{2n}}{3}.
\]

Since
\[
\sum_{i=1}^{n-1}a_i^2
\leq (n-1)a_{n-1}^2,
\]
we obtain
\[
(n-1)a_{n-1}^2
\geq \left(1-o(1)\right)\frac{2^{2n}}{3}.
\]

Consequently,
\[
{
a_{n-1}\geq
\left(\frac{1}{\sqrt{3}}-o(1)\right)\frac{2^n}{\sqrt{n-1}}=\left(\frac{2}{\sqrt{3}}-o(1)\right)\frac{2^{n-1}}{\sqrt{n-1}}
},
\]
which is better than the lower bound obtained directly, i.e., \(a_{n-1} \geq \left(\sqrt{\frac{2}{\pi}}-o(1)\right)
\frac{2^{n-1}}{\sqrt {n-1}}.\)

\bibliographystyle{plain}
\bibliography{biblio}

\end{document}